\documentclass[11pt,leqeq]{article}
\usepackage{amssymb,amsthm}
\usepackage{amsmath}
\usepackage{amscd}
\usepackage{xy}
\xyoption{all}
\usepackage[dvips]{graphicx}
\newtheorem{thm}{Theorem}
\newtheorem{prop}[thm]{Proposition} 

\newtheorem{defi}[thm]{Definition}

\DeclareFontFamily{OMS}{rsfs}{\skewchar\font'60}
\DeclareFontShape{OMS}{rsfs}{m}{n}{<-5>rsfs5 <5-7>rsfs7 <7->rsfs10
}{} \DeclareSymbolFont{rsfs}{OMS}{rsfs}{m}{n}
\DeclareSymbolFontAlphabet{\scr}{rsfs}

\newcommand{\des}{\displaystyle}

\newcommand{\im}{{\rm{im}}}

\newcommand{\Co}{{\rm Cob}}

\newcommand{\Mfu}{{\rm MFunc}}

\newcommand{\Ho}{{\rm H}}

\newcommand{\ve}{{\rm vect}}
\newcommand{\sn}{{\rm S_{n}}}

\newcommand{\HL}{{\rm HLQFT}}

\newcommand{\ho}{{\rm Hom}}

\renewcommand{\des}{\displaystyle}

\def\C{{\cal{C}}}

\def\d{{\cal{D}}}
\def\f{{\cal{F}}}

\def\h{{\cal{H}}}

\begin{document}
\title{Membrane Topology}
\author{Edmundo Castillo and
Rafael D\'\i az}
 \maketitle

\begin{abstract}
We construct membrane homology groups $\h(M)$ associated with each
compact connected oriented smooth manifold, and show that $\h(M)$
is matrix graded algebra.\\

 \end{abstract}

\section{Introduction}

In this paper we continue our research of homological quantum
field theories \HL \ initiated in \cite{Cas1} focusing our
attention in the two dimensional situation.  Our definition of
 \HL \ is based on several sources. The leading actor is the category $\Co_d$ of $d$ dimensional cobordisms. Objects in $\Co_d$
are boundaryless compact oriented smooth manifolds. Morphisms in
$\Co_d,$ called cobordisms, from $P$ to $Q$ are diffeomorphisms
classes of compact oriented $d-1$ dimensional manifolds with
boundary $M$ together with a diffeomorphism from $ (P^{-}
\sqcup Q) \times [0,1)$ onto an open neighborhood of
$\partial M$. Composition is given by gluing of cobordisms along
their boundaries.\\

Gradually it has become clear that the geometric background  for
quantum fields are monoidal representations of the category of
cobordisms, i.e., monoidal functors $\f \colon\Co_d
\rightarrow \ve.$\footnote{There are additional
constrains for a realistic quantum field theory other than those
imposed by the fact that they yield monoidal representations of
$\Co_d $. } Different types of field theories correspond with
different data on objects and morphisms in the cobordisms
category. For example in full quantum field theories objects and
cobordisms are given Riemmanian or Lorentzian metrics. Similarly,
in conformal field theory \cite{Se} cobordisms are endowed with
Riemmanian metrics defined up to conformal equivalence.\\

Topological quantum field theories, an important tool in modern
algebraic topology, are objects in the category of monoidal
functors $\f \colon\Co_d \rightarrow \ve$ introduced by Atiyah
\cite{MA}. Turaev in \cite{TU1} and \cite{TU2} introduced the notion of homotopical
quantum field theory. Fix a compact connected smooth manifold $M$.
Objects in the category $\mathrm{H}\Co_d^{M}$ of homotopically
extended cobordisms in $M$ are $d-1$ dimensional smooth compact
manifolds $N$ together with a homotopy class of maps $f:N
\rightarrow M$. A morphism in $\mathrm{H}\Co_d^{M}$ from $N_0$ to $N_1$ is
a cobordism $P$ connecting $N_0$ and $N_1$ together with a
homotopy class of  maps $g:P
\rightarrow M$ such that its restriction to the boundary
gives the homotopy classes associated with the boundary maps. A
homotopical quantum field theory is a monoidal functor
$\mathrm{H}\Co_d^{M} \rightarrow \ve$.\\

In this note we work within the context of homological quantum
field theory $\mathrm{HLQFT}$, i.e., monoidal representatios of
$\Co_d^{M}$ the category of homological extended cobordisms in
$M$. Objects in $\Co_d^{M}$ are $d-1$ dimensional manifolds $N$
together with a map sending each boundary component of $N$ into an
oriented embedded submanifold of $M$. Morphisms are cobordisms
together with an homology class of maps \footnote{Each map should
constant on a neighborhood of each boundary component, and mapping
each boundary component into its associated embedded submanifold
of $M$.}  from the cobordism into $M$. Composition of morphisms in
$\Co_d^{M}$ combines the usual composition of cobordisms with a
new sort of techniques introduced in the context of string topology.\\

Chas and Sullivan in their seminal paper \cite{SCh} showed that
the homology of the space of free loops on compact connected
oriented manifolds, after a degree shift, carries the structure of
a Batalin-Vilkovisky algebra. In order to define the loop product,
the $\triangle$ operator, and the bracket a sort of intersection
product on the homology of infinite dimensional manifolds is
required. Cohen and Jones, see
\cite{CG}, \cite{CoVo}, \cite{cj} and \cite{cjj} , showed that the
construction of the intersection product for certain infinite
dimensional manifolds reduces to the construction of the umkehr
map $F_{!}:H(\mathcal{N})
\longrightarrow H(\mathcal{M}),$ for  $F$ a smooth map between infinite dimensional
manifolds $\mathcal{N}$ and $\mathcal{M}$.  They proved that the
umkehr map exists if $F$ is a regular embedding of finite
codimension, which occurs if $F$ fits into a pullback square
diagram
\[\xymatrix @R=.3in  @C=.4in
{\mathcal{N} \ar[r]^{F} \ar[d]_{p} &  \mathcal{M}\ar[d]^{q}
\\ N \ar[r]^{f} & M  } \]
where $N$ and $M$ are finite dimensional manifolds, $f$ is an
embedding and $q$ is a fiber bundle. In this paper we use this
fundamental fact repeatedly, and refer to it as the Cohen and
Jones technique.\\

Our goal in this paper is to study homological quantum field
theories, i.e., monoidal functors $\Co_d^{M} \rightarrow
\ve$ using the Cohen and Jones technique. In Section
\ref{s0} we formally defined and give an example of \HL \ . In Section \ref{s1} we describe one dimensional
\HL \ using Cohen and Jones technique and show that one can
construct examples from connections on fiber bundles. In Section
\ref{sec3} we associate a two dimensional
\HL \ with each $B$-field defined on a connected oriented smooth
manifold. In Section \ref{sec4} we construct membrane homology
groups $\h(M)$ associated with each compact connected oriented
smooth manifold and show that it is matrix graded algebra equipped
with a natural representation.

\section{Homological quantum field theory}\label{s0}

Homology groups of space $M$ are denoted by $H(M)$. Let
$\Ho(M)=H(M)[d]$ be the homology of $M$ with degrees shifted by $d
\in \mathbb{N}$, i.e., $\Ho_i(M)=H_{i+d}(M)$. Let $D(M)$ be the set of
connected oriented embedded submanifolds of $M.$ The empty set is
assumed to be a $d$-dimensional manifold for $d\in
\mathbb{N}.$ \\

Objects in the category $\Co_{d}^{M}$ of homologically extended
cobordisms are triples $(N,f,<)$ such that $N$ is a compact
oriented manifold of dimension $d-1$,  $f\colon
\pi_{0}(N) \to D(M)$  is a map, and $<$ is a linear ordering on $\pi_{0}(N).$
We use the notation $\overline{f}=
\prod_{c\in \pi_{0}(N)}f(c)$. \\

For objects $(N_{0},f_{0},<_{0})$ and  $(N_{1}, \ f_{1},<_{1})$ in
$\Co_{d}^{M}$ we set \[ \Co_{d}^{M}((N_{0},f_{0},<_{0}), \
(N_{1},f_{1},<_{1}))=
\overline{\Co_{d}^{M}}((N_{0},f_{0},<_{0}), \ (N_{1},f_{1},<_{1}))
\diagup \backsim,\] where $\overline{\Co_{d}^{M}}((N_{0},f_{0},<_{0}), \
(N_{1},f_{1},<_{1}))$
is the set of triples $(P,\alpha,c)$ such that\\
\begin{itemize}
\item{$P$ is a compact oriented $d$-manifold with boundary.}

\item{ $\alpha \colon N_{0}\bigsqcup N_{1}\times [0,1)\to
\im (\alpha)\subseteq P$ is a diffeomorphism  such that $\alpha|_{N_{0}}$ reverses orientation and
$\alpha|_{N_{1}}$ preserves orientation.}

\item{ $c \in \Ho(M^{P}_{f_{0},f_{1}})=H(M^{P}_{f_{0},f_{1}})[\dim(\overline{f_{1}})]$
, where $M^{P}_{f_{0},f_{1}}$ denotes the space of smooth maps
$g\colon P\to M$ such that $g$ is constant on a neighborhood of
each connected component of  $\partial P.$ }
\end{itemize}

Triples $(P,\alpha,\xi )$ and $(P',\alpha', \xi' )$ in
$\overline{\Co_{d}^{M}}((N_{0},f_{0},<_{0}), \
(N_{1},f_{1},<_{1}))$  are  $\sim$ equivalent if there exists an
orientation preserving diffeomorphism $\varphi\colon P_{1} \to
P_{2}$ such that $\varphi \circ
\alpha=\alpha'$ and $\varphi_{\star}(\xi)=\xi'.$

\begin{thm}
 {\em $(\Co_{n}^{M}, \ \sqcup, \ \emptyset)$ is a monoidal category
with product $\sqcup$ and unit $\emptyset$. }
\end{thm}

\begin{proof}
Assume we are given morphisms $(P,\alpha,c)\in
 \Co_{d}^{M}
((N_{0},f_{0},<_{0}), (N_{1}, f_{1}, <_{1}))$ and $(Q,\beta,d)\in
 \Co_{d}^{M}
((N_{1},f_{1},<_{1}), (N_{2}, f_{2}, <_{2})).$ The composition
morphism $$(P,\alpha,c)\circ (Q,\beta,d)\in
\Co_{d}^{M}((N_{0},f_{0},<_{0}),
(N_{2},f_{2},<_{2}))$$  is the triple $(P \circ Q,\alpha \circ
\beta,c \circ d)$ where
$ P \circ Q = P \des\bigsqcup_{N_{1}} Q,$  $\alpha
\circ
\beta = \alpha \mid_{N_{2}}\bigsqcup
\beta \mid_{N_{0}},$ and $c \circ d$ is constructed from
the pull back diagram
\[\xymatrix @R=.2in  @C=.5in
{M^{P}_{f_{0},f_{1}}\times_{\overline{f_{1}}} M^{Q}_{f_{1},f_{2}}
\ar[r]^{d} \ar[d]_{e} & M^{P}_{f_{0},f_{1}}\times
M^{Q}_{f_{1},f_{2}}\ar[d]^{e_{t}\times e_{s}}
\\ \overline{f_{1}} \ar[r]^{\Delta} & \overline{f_{1}} \times \overline{f_{1}}} \]
as the composition of maps \[\Ho( M^{P}_{f_{0},f_{1}})\otimes \Ho(
M^{Q}_{f_{1},f_{2}}) \stackrel{d!}\longrightarrow
\Ho(M^{P}_{f_{0},f_{1}}\times_{\overline{f_{1}}}  M^{Q}_{f_{1},f_{2}}) \stackrel{i_{\star}}\longrightarrow
\Ho(M^{P\sqcup_{N_{1}} Q}_{f_{0},f_{2}})\] where $i: M^{P}_{f_{0},f_{1}}\times_{\overline{f_{1}}}
M^{Q}_{f_{1},f_{2}} \longrightarrow M^{P\sqcup_{N_{1}}
Q}_{f_{0},f_{2}}$ sends a pair $(x,y)$ into the map
$i(x,y):P\sqcup_{N_{1}} Q \longrightarrow M$ whose restriction to
$P$ is $x$ and whose restriction to $Q$ is $y$. Associativity is
proved as in the case of string topology \cite{cj} .The identity
morphism is $(N \times [0,1] ,\alpha,1_N)\in \Co_{d}^{M}( (N,f,<),
(N, f, <))$, where $1_N$ is defined as follows: consider the map
$c:N \longrightarrow  M_{f,f}^{N \times [0,1]}$ sending $n \in N$
to the map constantly equal to $n$, then $1_N = c_*([N])$.

\end{proof}

\begin{defi}
 {\em $(\Co_{n,r}^{M}, \ \sqcup)$ is the full monoidal subcategory
of $\Co_{n}^{M}$ without unit. }
\end{defi}

Given monoidal categories $\C$ and $\d$ we let $\Mfu(\C,\d)$ be
the category of monoidal functors from $\C$ to $\d.$

\begin{defi}
{\em  The category of $d$ dimensional homological quantum field
 theories is given by $\HL_{d}(M)=\Mfu(\Co_{n}^{M},\ve)$.  The category of the $d$
dimensional restricted homological quantum field theories is
$\HL_{d,r} (M)=\Mfu(\Co_{d,r}^{M},\ve)$.}
\end{defi}

Let us construct an example of restricted homological quantum
field theory.

\begin{thm} {\em The map $\Ho \colon
\Co_{d,r}^{M} \to
\mbox{vect}$ given on objects by  $\Ho(N,f,<)= \Ho(\overline{f})$
defines a restricted homological quantum field theory.}
\end{thm}
\begin{proof}

Fix $P$ a cobordism between $N_0$ and $N_1$. We need a map $
\Ho(M_{f_0,f_1}^P) \rightarrow \ho(\Ho(\overline{f_0}),\Ho(\overline{f_1})),$
or equivalently, an adjoint map $\Ho(\overline{f_{0}})\otimes
\Ho( M^{P}_{f_{0},f_{1}}) \longrightarrow \Ho(\overline{f_{1}}).$ The pullback diagram
\[\xymatrix @R=.2in  @C=.5in
{\overline{f_{0}}\times_{\overline{f_{0}}} M^{P}_{f_{0},f_{1}}
\ar[r]^{d} \ar[d]_{e} &
\overline{f_{0}}\times M^{P}_{f_{0},f_{1}} \ar[d]^{I \times e_{s}}
\\ \overline{f_{0}} \ar[r]^{\Delta} & \overline{f_{0}} \times \overline{f_{0}}} \]
induces the desired map through the compositions \[
\Ho(\overline{f_{0}})\otimes \Ho( M^{P}_{f_{0},f_{1}})
\stackrel{d!}\longrightarrow
\Ho(\overline{f_{0}}\times_{\overline{f_{0}}} M^{P}_{f_{0},f_{1}}) \stackrel{t_{\star}}\longrightarrow
\Ho(\overline{f_{1}}).\] Units and associativity are constructed
as in the previous theorem.
\end{proof}

\section{One dimensional homological quantum field theory }\label{s1}

In this section we study  \HL \ in dimension one using the Cohen
and Jones technique.  For a manifold $N$ we let $\pi_{0}(N)$ be
the set of connected components of $N$, and we set
$\overline{N}=\prod_{c\in
\pi_{0}(N)}c.$ Objects in open string category \cite{DS} are embedded submanifolds of
$M$. For $N_0$ and $N_1$ embedded submanifolds of $M$ the space of
morphisms is $\Ho(M^{I}_{N_0,N_1})$, where $M_{N_0, N_1}^{I}$ be
the space of smooth maps $x:I
\longrightarrow M$  constant on  neighborhoods of
$0$ and $1$, respectively. Let $\Ho(M_{N_0, N_1}^{I})=H(M_{N_0,
N_1}^{I})[\dim(N_1)]$. Composition of morphisms is defined as
follows. We have a pullback diagram
\[\xymatrix @R=.2in  @C=.5in
{M_{N_0, N_1}^{I}\times_{N_1} M_{N_1, N_2}^{I} \ar[r]^{d}
\ar[d]_{e} & M_{N_0, N_1}^{I}\times M_{N_1,
N_2}^{I}\ar[d]^{e\times e}
\\ N_1 \ar[r]^{\Delta} & N_1 \times N_1  } \]
and a map $i:M_{N_0, N_1}^{I}\times_{N_1} M_{N_1, N_2}^{I}
\longrightarrow M_{N_0, N_2}^{I}$ sending a pair $(x,y)$ to the
path that runs trough $x$ in half the time and then trough $y$ in
the the other half. Consider the following map $\bullet:
\Ho(M_{N_0, N_1}^{I})\otimes \Ho(M_{N_1, N_2}^{I}) \longrightarrow
\Ho(M_{N_0, N_2}^{I}) $ given through compositions
\[\Ho(M_{N_0, N_1}^{I})\otimes \Ho(M_{N_1, N_2}^{I})
\stackrel{d!}\longrightarrow
\Ho(M_{N_0, N_1}^{I}\times_{\overline{N_{1}}} M_{N_1, N_2}^{I}) \stackrel{i_{\star}}\longrightarrow
\Ho(M_{N_0, N_2}^{I}).\]

Let us now consider \HL \ in dimension one.  An object $f$ in
$\Co^{M}_{1,r}$ is a map $f
\colon [n]\to D(M)$ where For $n\in
\mathbb{N}^+$ we set $[n]=\{1,\cdots ,n\}.$  We use the notation
$\overline{f}= \prod_{i\in[n]}f(i).$  The space of morphisms in
$\Co_{1,r}$ from $f$ to $g$ is
\[\Co^{M}_{1,r}(f,g)= \bigoplus_{\sigma\in \sn} \bigotimes_{i=1}^{n}
\Ho(P(f(i),g(\sigma(i)))).\] Composition of morphisms in $\Co_{1,r}^{M}$ is given by
\[\xymatrix @R=.15in{\Co^{M}_{1,r}(f,g)\otimes \Co^{M}_{1,r}(g,h) \ar[d]^{=}\\
\bigoplus_{\sigma, \tau\in \sn} \bigotimes_{i=1}^{n}
\Ho(P(f(i),g(\sigma(i))))\otimes \bigotimes_{j=1}^{n}
\Ho(P(g(j),h(\tau(j)))) \ar[d]^{s}\\
 \bigoplus_{\sigma, \tau\in \sn} \bigotimes_{i=1}^{n}
\Ho(P(f(i),g(\sigma(i))))\otimes
\Ho(P(g(\sigma(i)),h(\tau(\sigma(i)))))\ar[d]^{\bullet}\\
\bigoplus_{\sigma, \tau\in \sn} \bigotimes_{i=1}^{n}
\Ho(P(f(i),h(\tau(\sigma(i)))))\ar[d]^{=}\\
\bigoplus_{\rho\in \sn} \bigotimes_{i=1}^{n}
\Ho(P(f(i),h(\rho(i))))\ar[d]^{=}\\
\Co^{M}_{1,r}(f,h)}\] the map $s$ permutes order in the tensor products. \\

Let $G$ a compact Lie group and $\pi\colon P\to M$ a principal $G$
bundle over $M$. Let $\scr{A}_{P}$ denote the space of connections
on $P$ and $\Lambda
\in \scr{A}_{P}$. If $\gamma \colon I \to M$ is a
smooth curve on $M$ and $x\in P$ is such that $\pi(x)=\gamma(0),$
then we let $P_{\Lambda}(\gamma,x)$ be $\widehat{\gamma}(1)$ where
$\widehat{\gamma}$ is the horizontal lift of $\gamma$ with respect
to $\Lambda$ such that $\widehat{\gamma}(0)=x.$ Our next goal is
to prove the following result.

\begin{prop}\label{sl}
{\em There is a natural map $\Ho\colon \scr{A}_{P} \to
\HL_{1,r}(M).$}
\end{prop}

For each connection $\Lambda \in \scr{A}_{P}$ we construct a
functor $\Ho_{\Lambda}\colon
\Co_{1,r}^{M}\to \ve$ given on an object $f$ by
$$\Ho_{\Lambda}(f)=\Ho({P\mid_{f}})=H({P\mid_{f}})[\dim(\overline{f})],$$
where $P\mid_{f(i)}$ denotes the restriction of $P$ to
$f(i)\subseteq M$ and ${P\mid_{f}}=\des\prod_{i\in[n]}
P\mid_{f(i)}$. Proposition \ref{sl} follows from the next result.

\begin{thm}\label{t6}
{\em The map $\Ho_{\Lambda}\colon \Co_{1,r}^{M}\to\ve$ sending $f$
into $\Ho_{\Lambda}(f)$ defines an one dimensional restricted
homological quantum field theory.}
\end{thm}

\begin{proof}

In order to define  $\Ho_{\Lambda}\colon  \Co_{1,r}^{M}(f,g)
\longrightarrow  \ho(\Ho({P\mid_{f}}),
\Ho({P\mid_{g}}))$ we construct the
adjoint map $\Ho({P\mid_{f}})\otimes \Co^{M}_{1,r}(f,g)
\longrightarrow \Ho({P\mid_{g}}).$ The pullback diagram
\[\xymatrix @R=.2in  @C=.5in
{\bigsqcup_{\sigma\in \sn}\prod_{i=1}^{n}P\mid_{f(i)}
\times_{f(i)}P(f(i),g(\sigma(i))) \ar[d]_{e}
\ar[r]^{d} & \bigsqcup_{\sigma\in \sn}\prod_{i=1}^{n}P\mid_{f(i)}
\times P(f(i),g(\sigma(i))) \ar[d]^-{e\times e}\\
\overline{f}  \ar[r]_{\triangle} & \overline{f}\times \overline{f}}\]
together with the map \[i\colon  \bigsqcup_{\sigma\in
\sn}\prod_{i=1}^{n}P\mid_{f(i)}
\times_{f(i)}P(f(i),g(\sigma(i)))\stackrel{P_{\Lambda}}\longrightarrow
 \bigsqcup_{\sigma\in \sn}\prod_{i=1}^{n}P\mid_{g(\sigma(i))}\stackrel{s}\longrightarrow
\prod_{i=1}^{n}P\mid_{g(i)}\] allow us to define the desired map through the compositions

\[\xymatrix @R=.15in {\Ho({P\mid_{f}})\otimes \Co^{M}_{1,r}(f,g) \ar[d]^{k}\\
\bigoplus_{\sigma, \in \sn} \bigotimes_{i=1}^{n}
\Ho(P\mid_{f_{i}}\times_{f(i)} P(f(i),g(\sigma(i)))) \ar[d]^{d!}\\
 \bigoplus_{\sigma, \in \sn} \bigotimes_{i=1}^{n}
\Ho(P\mid_{f_{i}}\times P(f(i),g(\sigma(i))))\ar[d]^{i_{\star}}\\
\Ho({P\mid_{g}})
}\]
\end{proof}

There is a remarkable analogy between \HL \ in dimension one and
the algebra of matrices. In \cite{Cas2} this analogy is studied,
and several well-known constructions for matrices are generalized
to the homological context, among them the notion of Schur
algebras. Representations of homological Schur algebras are deeply
related with one dimensional quantum field theories. Further
examples of \HL \ in dimension one are considered in \cite{Cas3}.

\section{Two dimensional homological quantum field theory }\label{sec3}

Let $M$ be a compact oriented smooth manifold. According to Segal
\cite{Se1} a $B$ field, also known as a gerbe
with connection, is a complex Hermitian line bundle $L$ on
$M^{S^1}$, loops in $M$, equipped with a string connection. A
string connection is a rule that assigns to each surface with
boundaries $\Sigma$ and each map $y\colon \Sigma
\rightarrow M$ a transport operator
$B_{y} \colon L_{\partial(\Sigma)_{-}} \longrightarrow
L_{\partial(\Sigma)_{+}},$ where the extensions of $L$ to
$(M^{S^1})^{n}$ is defined by the rule $L_{(x_1,...,x_n)}=L_{x_1}
\otimes ... \otimes L_{x_{n}}.$ The assignment $y \rightarrow
B_y$ is assumed to have the following properties

\begin{itemize}
\item{It is a continuous map taking values in unitary operators. Therefore we have induced maps
$B_{y}\colon L^{1}_{\partial(\Sigma)_{-}} \rightarrow
L^{1}_{\partial(\Sigma)_{+}}$ between the corresponding circle
bundles.}

\item{It is transitive with respect to the gluing of surfaces.}

\item{It is parametrization invariant.}
\end{itemize}

Our next goal is to prove the following
\begin{prop}
{\em There is a natural map from $B\mbox{ fields on } M$ to
$\HL_{2,r}(M).$}
\end{prop}

Thus for each $B$ field we need a functor $\Ho_{B}\colon
\Co^{M}_{2,r} \longrightarrow \ve$. An object in $\Co^{M}_{2,r}$ is an
disjoint union of $n$ ordered circles together with a map $f \colon [n]
\longrightarrow D(M)$.  The functor $H_B$ is defined by the rule $\Ho_{B}(f)=
\Ho({L}_{f}^{1})= H({L}_{f}^{1})[\dim(\overline{f})]$, where $\Ho({L}_{f}^{1})=
\Ho(L^{1}\mid_{f(1)\times ... \times f(n)}).$ The notation
$L^{1}\mid_{f(1)\times ... \times f(n)}$ makes sense since
$f(1)\times ... \times f(n) \subseteq M
\times ... \times M \subseteq M^{S^1} \times ... \times M^{S^1}.$

\begin{thm}
{\em The map $\Ho_{B}\colon \Co_{2}^{M}\to\ve$ sending $f$ into
$\Ho(L_{f}^1)$ defines a two dimensional restricted homological
quantum field theory.}
\end{thm}

\begin{proof}
Consider the pullback diagram
\[\xymatrix @R=.2in  @C=.5in
{{L}_{f}^{1}\times_{\overline{f}}M_{_{f},_{g}}^{\Sigma}
 \ar[d]_{e}
\ar[r]^{d} & {L}_{f}^{1}
\times M_{_{f},_{g}}^{\Sigma} \ar[d]^-{\pi \times e}\\
\overline{f}  \ar[r]_{\triangle} & \overline{f}\times
\overline{f}}\] For $f$ and $g$  objects in $\Co^{M}_{2,r}$, the
$B$
field induces a map ${L}_{f_{0}}^{1}
\times_{\overline{f}}M_{_{f},_{g}}^{\Sigma}
\stackrel{e_B}\longrightarrow {L}_{g}^{1}.$ We need a map
$\Ho(M_{_{f},_{g}}^{\Sigma}) \longrightarrow
Hom(\Ho({L}_{f}^{1}),\Ho({L}_{g}^{1}))$, its adjoint map is given
by the compositions
\[ \Ho({L}_{f}^1)\otimes \Ho(M_{_{f},_{g}}^{\Sigma})
\stackrel{d!}\longrightarrow
\Ho({L}_{f}^ 1 \times_{\overline{f}} M_{_{f},_{g}}^{\Sigma}) \stackrel{e_{B*}}\longrightarrow
\Ho(L_{g}^1). \]
\end{proof}

\section{Membrane topology} \label{sec4}

Let us take a closer look at objects in the category
$\Co^{M}_{2,r}$. We focus our attention on objects $f
\colon [n] \longrightarrow D(M)$ such that $f$ is constantly equal
to $M$, and so objects are just integers. A morphism from $n$ to
$m$ is a homology class of $M^{\Sigma}$ where $\Sigma$ is a
compact oriented surface with $n$ incoming boundary components and
$m$ outgoing boundary components. We further restrict our
attention to connected surfaces $\Sigma$.\\

\begin{defi}
{\em For integers $n,m \geq 1,$ let $\Sigma_{n,g}^{m}$ be a
Riemann surface of genus $g$ with $n$ incoming numbered marked
points and $m$ outgoing numbered marked points.}
\end{defi}

\begin{figure}[ht]
\begin{center}
\includegraphics[height=4cm]{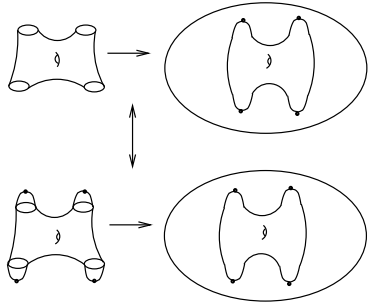}
\caption{ \ An element in $M^{\Sigma}$ and the corresponding element in $M^{\Sigma_{n}^{m}}$ . \label{bra1}}
\end{center}
\end{figure}

Let $M^{\Sigma_{n,g}^{m}}$ be the space of smooth maps $x\colon
\Sigma_{n,g}^{m}
\longrightarrow M$ constant in a neighborhood of each marked
point. If $\Sigma$ is a genus $g$ surface with $n$ incoming
boundaries and $m$ outgoing boundaries, then the spaces
$M^{\Sigma}$ and $M^{\Sigma_{n,g}^{m}}$ are homotopically
equivalent, see Figure \ref{bra1}, and therefore
$H(M^{\Sigma})=H(M^{\Sigma_{n,g}^{m}}).$

Let us introduce an algebraic notion.

\begin{defi}
{\em Algebra $(A,m)$ is a matrix graded if $A=
\bigoplus_{n,m=1}^{\infty} A_{n}^{m},$ $m\colon A_{n}^{m}\otimes A_{m}^{k}\to A_{n}^{k},$
and  $m\mid_{A_{n}^{m}\otimes A_{p}^{k}}=0$ if $p\neq m$.}
\end{defi}

\bigskip

We are ready to define  membrane homology groups.

\begin{defi} {\em Membrane homology of a compact oriented manifold $M$ is given by
$\mathcal{H}(M)=\bigoplus_{n,m=1}^{\infty}
\Ho_{n}^{m}(M),$ where $\Ho_{n}^{m}(M)=\bigoplus_{g=0}^{\infty} \Ho^{m}_{n,g}(M)$ and
$\Ho^{m}_{n,g}(M)=H(M^{\Sigma_{n,g}^{m}})[m\dim M].$}
\end{defi}

\begin{thm}
{\em $\h(M)$ is a matricially graded algebra.}
\end{thm}

\begin{proof}
There is a pullback square diagram
\[\xymatrix @R=.3in  @C=.5in
{M^{\Sigma_{n,g}^{m}}\times_{M^{m}}M^{\Sigma_{m,h}^{k}}
 \ar[d]_{e}
\ar[r]^{d} & M^{\Sigma_{n,g}^{m}}\times M^{\Sigma_{m,h}^{k}} \ar[d]^-{e_{t}\times e_{s}}\\
M^{m}  \ar[r]^{\triangle} & M^{m}\times M^{m}}\] and a natural map
$i\colon M^{\Sigma_{n,g}^{m}}\times_{M^{m}}M^{\Sigma_{m,h}^{k}}\to
M^{\Sigma_{n,g+h+m-1}^{k}}$  which is better explained by Figures
$\ref{cobra2}$,  $\ref{cobra1}$ and $\ref{cobra3}$  below, where a
pair $(x,y)$ in
$M^{\Sigma_{3,2}^{2}}\times_{M^{m}}M^{\Sigma_{2,1}^{1}}$ is shown
as well as the induced element $i(x,y) \in M^{\Sigma_{3,4}^{1}}.$

\begin{figure}[h]
\begin{minipage}[t]{0.5\linewidth}
\begin{center}
\includegraphics[height=4cm]{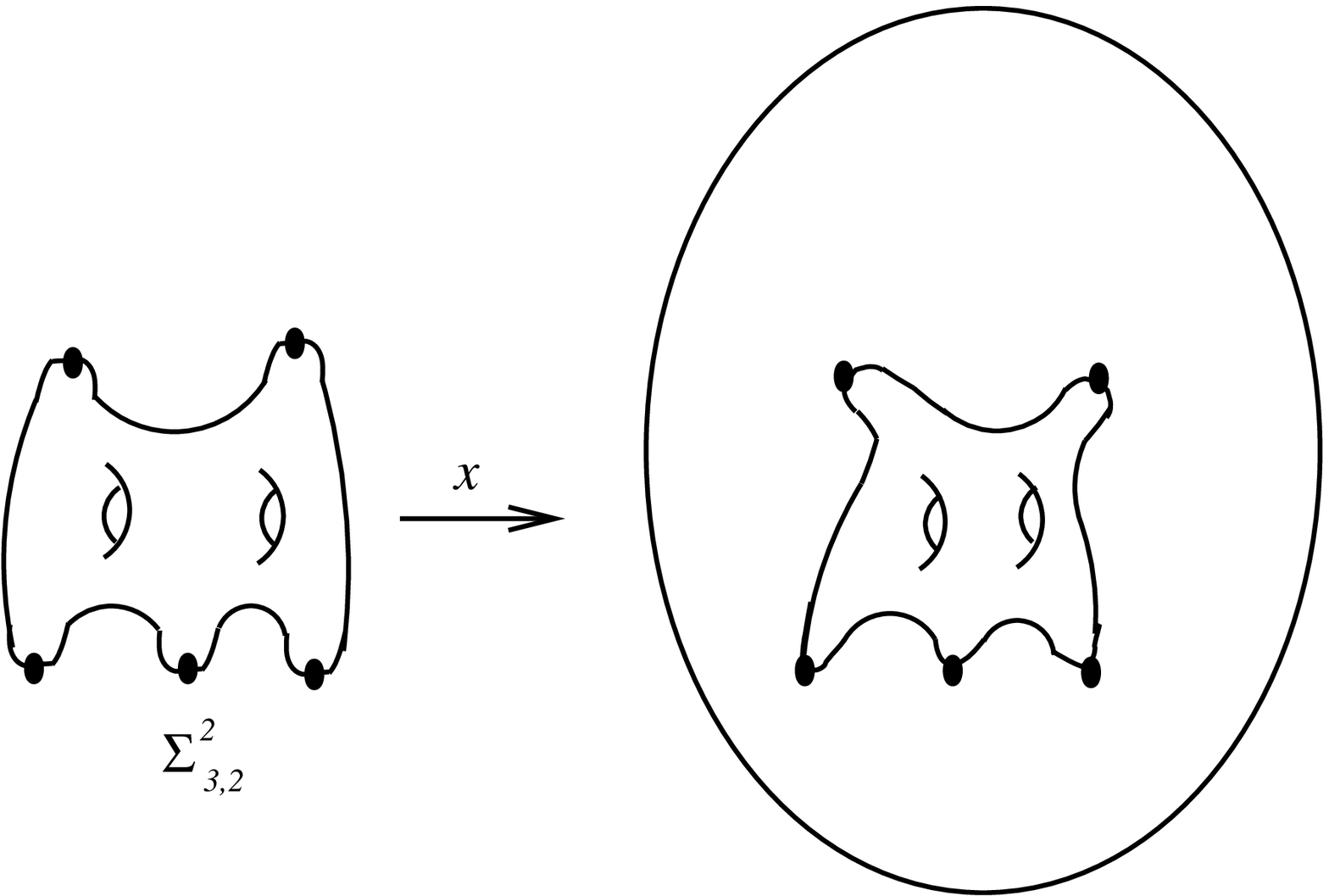} \caption{ \ Element in $M^{\Sigma_{3,2}^{2}}$.}\label{cobra2}
\end{center}
\end{minipage}
\begin{minipage}[t]{0.5\linewidth}
\begin{center}
\includegraphics[height=4cm]{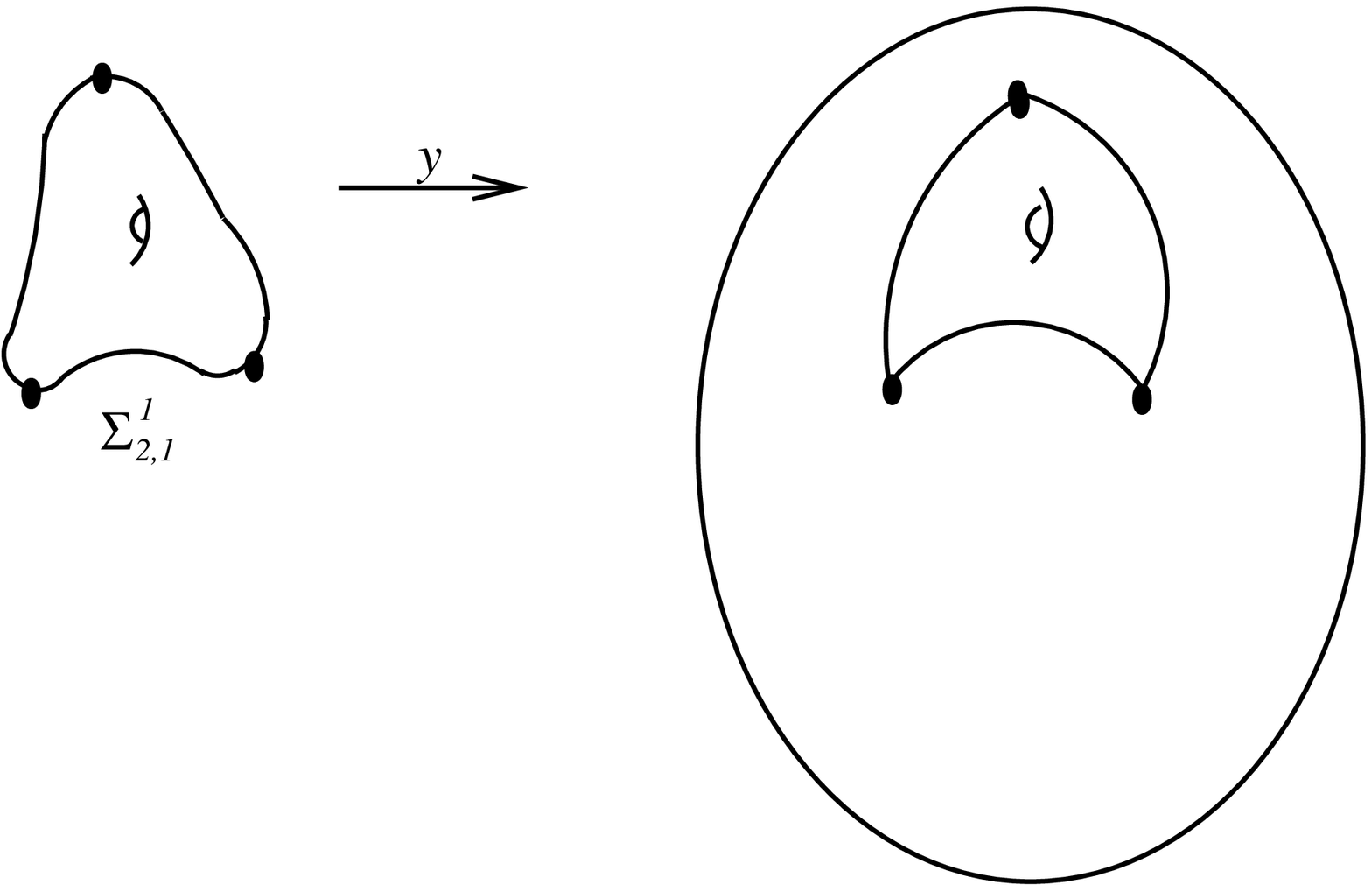} \caption{ \ Element in $M^{\Sigma_{2,1}^{1}}$.}\label{cobra1}
\end{center}
\end{minipage}
\end{figure}

\begin{figure}[ht]
\begin{center}
\includegraphics[height=4cm]{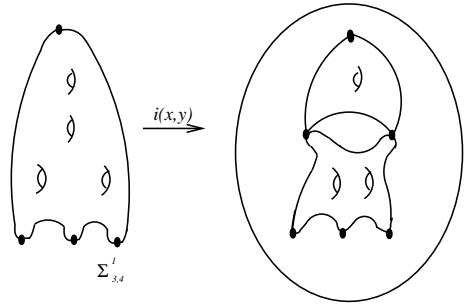}
\caption{ \ Element $i(x,y) \in M^{\Sigma_{3,4}^{1}}$. \label{cobra3}}
\end{center}
\end{figure}

From the pullback diagram and the map $i$ above, we define a
product on $\mathcal{H}(M)$ via the composition of maps
\[ \Ho(M^{\Sigma_{n,g}^{m}})\otimes \Ho(M^{\Sigma_{m,h}^{k}})
\stackrel{d!}\longrightarrow
\Ho(M^{\Sigma_{n,g}^{m}}\times M^{\Sigma_{m,h}^{k}}) \stackrel{i_{\star}}\longrightarrow
\Ho(M^{\Sigma_{n,g+h+m-1}^{k}}) \]
\end{proof}

Next we show that membrane homology comes equipped with a natural
representation. For a vector space $V$ we let $T_{+}(V)=\des
\bigoplus_{n=1}^{\infty}V^{\otimes n}.$

\begin{thm}
{\em $T_{+}(\Ho(M))$ is a representation of $\h(M)$.}
\end{thm}

\begin{proof}

The pullback diagram

\[\xymatrix @R=.2in  @C=.5in
{M^{n}\times_{M^{n}}M^{\Sigma_{n,g}^{m}} \ar[d]_{e}
\ar[r]^{d} & M^{n}\times M^{\Sigma_{n,g}^{m}} \ar[d]^-{I \times e}\\
M^{n}  \ar[r]^{\triangle} & M^{n}\times M^{n}}\] and the map
$i\colon M^{n}\times_{M^{n}}M^{\Sigma_{n,g}^{m}}\to M^{m}$, induce
an action of $\h(M)$ on $T_{+}(\Ho(M))$ via the composition of
maps

\[ \Ho(M)^{\otimes n}\otimes \Ho(M^{\Sigma_{n,g}^{m}})
\stackrel{d!}\longrightarrow
\Ho(M^{m}\times_{M^{n}} M^{\Sigma_{n,g}^{m}}) \stackrel{i_{\star}}\longrightarrow
\Ho(M)^{\otimes m} \]

\end{proof}
As we have seen membrane topology is an interesting algebraic
structure associated with any oriented manifold. It would be
interesting to compute it explicitly for familiar spaces, and also
to study its relations with other types of two dimensional field
theories, such as topological conformal field theories in the
sense of \cite{KVS} and
\cite{Kimm}.

\subsection*{Acknowledgment} Thanks to Jaime Camacaro, Takashi Kimura, Lorenzo Leal, Eddy Pariguan,
Bernardo Uribe, Raymundo Popper and Arturo Reyes.

\bibliographystyle{amsplain}
\bibliography{SurfacesTopology}

\noindent ragadiaz@gmail.com, \ \ ecastill@euler.ciens.ucv.ve\\
\noindent Universidad Central de Venezuela, Caracas, Venezuela.\\

\end{document}